\documentclass{article}
\usepackage[utf8]{inputenc}

\usepackage{amsmath}
\usepackage{amsthm}
\usepackage{amssymb}
\usepackage{wasysym}
\usepackage[T1]{fontenc}
\usepackage{graphicx}
\usepackage{caption}
\usepackage[english]{babel}
\usepackage[pdfencoding=auto]{hyperref}
\newtheorem{theorem}{Theorem}

\newtheorem{lemma}{Lemma}

\newtheorem{assumption}{Assumption}
\newtheorem{definition}{Definition}
\newtheorem{remark}{Remark}
\usepackage{color} 
\usepackage{url}
\everymath{\displaystyle}

\title{\Large \bf
On systematic criteria for the global stability of nonlinear systems via the Koopman operator framework
}

\author{Christian Mugisho Zagabe and Alexandre Mauroy
\thanks{C.M. Zagabe is with Department of mathematics,
        University of Namur, 5000 Namur, Belgium
        {\tt\small christian.mugisho@unamur.be}}%
\thanks{A. Mauroy is with Department of mathematics,
        University of Namur, 5000 Namur, Belgium
        {\tt\small alexandre.mauroy@unamur.be}}%
}
\date{}

\begin{document}

\maketitle


\begin{abstract}
We present novel sufficient conditions for the global stability of equilibria in the case of nonlinear dynamics with analytic vector fields. These conditions provide stability criteria that are directly expressed in terms of the Taylor expansion coefficients of the vector field (e.g. in terms of first order coefficients, maximal coefficient, sum of coefficients). Our main assumptions is that the flow be holomorphic, and the linearized system be locally exponentially stable and diagonalizable. These results are based on the properties of the Koopman operator defined on the Hardy space on the polydisc.

\end{abstract}

\section{Introduction}

In dynamical systems theory, characterizing global stability remains a challenge. The existence of a Lyapunov function guarantees global stability due to Lyapunov's second method, but there are very few general constructive methods. For a linear system, on the other hand, the existence of a quadratic Lyapunov function is both a necessary and sufficient condition for global stability. 
In this context, the Koopman operator approach provides a ``global linearization" of nonlinear dynamics (see e.g., \cite{MBRMIM, AMIMYS}), which is amenable to global stability analysis through linear methods \cite{AMIMYS}. 
For instance, specific Koopman eigenfunctions were used in \cite{AMIM1} to obtain necessary and sufficient conditions for global stability of hyperbolic attractors, a result which mirrors well-known spectral stability results for linear systems. Moreover, a numerical method was proposed in \cite{MauroySootlaMezic_stability} to compute Lyapunov functions  from a finite dimensional approximation of the Koopman operator. 

The present work follows the same path as the above-mentioned results based on the Koopman operator approach. However, it does not use Koopman eigenfunctions, which are usually unknown and have to be computed numerically, nor does it rely on possibly inaccurate approximations of the operator. Instead, our results provide sufficient conditions for global stability of equilibria associated with holomorphic flows, which can be directly verified with the system vector field.
Under mild assumptions on the linearized dynamics (i.e. exponential stability, diagonalizable linear system), the specific case of polynomial vector fields is considered, along with more general analytic vector fields.
Our theoretical findings are built upon our previous work based on the properties of the Koopman generator defined in the Hardy space of the polydisc \cite{CMZAM2}. But in contrast to previous work, they do not focus on switched nonlinear systems and provide stability conditions which are less conservative thanks to the use of re-scaled Hardy spaces. Moreover, the obtained criteria are more readily applicable since they are expressed in terms of simple quantities directly computed from Taylor coefficients (e.g. first order coefficients, maximal coefficient, discounted sum of coefficients).

The remainder of the paper is structured as follows. In Section \ref{sec:prelim}, we provide a general introduction to the Koopman operator framework and some specific properties in the Hardy space on the polydisc. Our main results are presented in Section \ref{sec:mainresult} and illustrated in Section \ref{sec:illustra} with two examples. Section \ref{sec:concl} gives concluding remarks and perspectives. The proofs of our main results can be found in Appendix \ref{sec:appendix_proof1} and \ref{sec:appendix_proof2}.

\subsection*{Notations}
For multi-index notations  $\alpha = (\alpha_1 ,..., \alpha_n)\in \mathbb{N}^n$, we define $\vert\alpha\vert = \alpha_1+ \cdots + \alpha_n$ and $z^\alpha = z_1^{\alpha_1}\cdots z_n^{\alpha_n}$.
The complex conjugate and real part of a complex number $a$ is denoted by $\bar{a}$ and $\Re(a)$, respectively. 
The Jacobian matrix of the vector field $F$ at $z$ is given by $JF(z)$.
The complex polydisc centered at $0$ and of radius $\rho>0$ is defined by 
$$\mathbb{D}^n(\rho)=\left\lbrace z\in \mathbb{C}^n: | z_1|<\rho,\cdots, | z_n|<\rho  \right\rbrace$$
and $\partial \mathbb{D}^n(\rho)$ and $\left(\partial \mathbb{D}(\rho)\right)^n$ is its boundary and distinguished boundary respectively.
In particular, $\mathbb{D}^n$ denotes the unit polydisc (i.e. with $\rho=1$).

\section{Preliminaries}
\label{sec:prelim}

We consider a continuous-time dynamical system
\begin{equation}\label{eq:systnon}
\dot{z} = F(z), \quad z\in \mathbb{D}^n(\rho), 
\end{equation}
with $\rho>0$, where the vector field $F$ satisfies the following assumption.
\begin{assumption}\label{assump4}
The components $F_l,\, l=1,\cdots, n$, of the vector field $F$ (i) are holomorphic on the closed polydisc $\overline{ \mathbb{D}^n(\rho)}$, (ii) belong to the Hardy space $\mathbb{H}^2(\mathbb{D}^n(\rho))$ (defined in Section \ref{sec:hardy_space} below), and 
(iii) generate a flow $\varphi^t$ that is holomorphic on $\mathbb{D}^n(\rho)$ and maps $\mathbb{D}^n(\rho)$ to itself.
\end{assumption}
\smallskip
Necessary and sufficient conditions on the vector field to ensure that Assumption \ref{assump4}(iii) is satisfied are given in \cite{RYCZHZ} for the case of the unit polydisc.

Moreover, we will make the following additional standing assumption related to the type of dynamical behavior we investigate.
\begin{assumption}
\label{assump_hyperbolic}
The vector field $F$ admits on $\mathbb{D}^n(\rho)$ a unique hyperbolic equilibrium at the origin (without loss of generality), i.e. $F(0)=0$, and the eigenvalues $\tilde{\lambda}_j$ of the Jacobian matrix $JF(0)$ satisfy $\Re\{\tilde{\lambda}_j\}<0$.
\end{assumption}

In order to investigate the global stability properties of the above dynamical system, we will define the Koopman operator on a proper space adapted to the dynamics. Since we made the assumption of analyticity of vector fields and flows, it is natural to consider a space of analytic functions, and a prototypical choice is the Hardy space on the polydisc.

\subsection{Hardy space of the polydisc}
\label{sec:hardy_space}

The \emph{Hardy space} of holomorphic functions on the \emph{ polydisc} $\mathbb{D}^n(\rho)$ is the space
 $$\mathbb{H}^2(\mathbb{D}^n(\rho))=\left\lbrace f:\mathbb{D}^n(\rho)\rightarrow \mathbb{C}, \mbox{holomorphic}:\|f\|^2<\infty\right\rbrace, $$
where 
\[ \|f\|^2=\lim_{r\rightarrow 1^-} \int_{(\partial \mathbb{D}(\rho))^n}\vert f\left(r\omega\right)\vert^2dm_n(\omega)\]
 and $m_n$ is the normalized Lebesgue measure on $(\partial \mathbb{D}(\rho))^n$. 
The space is equipped with an inner product defined by 
$$\left\langle f,g\right\rangle=\int_{(\partial \mathbb{D}(\rho))^n}f\left(\omega\right) \bar g\left(\omega\right) dm_n(\omega),$$
so that the set of monomials $\left\lbrace z^\alpha: \alpha \in  \mathbb{N}^n\right\rbrace$ is a standard orthonormal basis on $\mathbb{H}^2(\mathbb{D}^n(\rho))$. In the sequel, the monomials will be denoted by $e_{k}(z)=z^{\alpha(k)}$, where the map $\alpha:\mathbb{N} \to \mathbb{N}^n$, $k\mapsto \alpha (k)$ refers to the lexicographic order\footnote{that is, $e_{k_1}<e_{k_2}$ if $|\alpha(k_1)|<|\alpha(k_2)|$, or if $|\alpha(k_1)|=|\alpha(k_2)|$ and $\alpha_j(k_1)>\alpha_j(k_2)$ for the smallest $j$ such that $\alpha_j(k_1)\neq\alpha_j(k_2)$}. 
 For $f$ and $g$ in $\mathbb{H}^2(\mathbb{D}^n(\rho)) $, with $f=\sum_{k \in \mathbb{N}} f_k e_k$ and $g=\sum_{k \in \mathbb{N}} g_k e_k$, the isomorphism 
$$\sum_{k \in \mathbb{N}} f_k e_k \mapsto (f_k)_{k \geq 0}$$  
 between $\mathbb{H}^2(\mathbb{D}^n(\rho))$ and the $l^2$-space allows to rewrite the norm and the inner product as
$$\Vert f\Vert^2=\sum_{k \in \mathbb{N}}\rho^{2|\alpha(k)|} \vert f_k \vert^2 \quad \mbox{and}\quad \left\langle f,g\right\rangle=\sum_{k \in \mathbb{N}} \rho^{2|\alpha(k)|} f_k \, \bar g_k.$$ 
By using the change of variables $\phi(z)=z'=z/\rho$ on $\mathbb{D}^n(\rho)$ , the map $f\mapsto f'=f\circ \phi$ defines an isometry between the two Hardy spaces  $\mathbb{H}^2(\mathbb{D}^n(\rho))$ and $\mathbb{H}^2(\mathbb{D}^n)$ where $\Vert f'\Vert^2_{\mathbb{H}^2(\mathbb{D}^n)}=\sum_{k \in \mathbb{N}} \vert f_k \vert^2 $. For more details on the Hardy space, we refer the reader to \cite{WR1,WR2,JHS}.

\subsection{Koopman operator on $\mathbb{H}^2(\mathbb{D}^n(\rho))$}

The \emph{Koopman operator} is defined here as the composition operator on $\mathbb{H}^2(\mathbb{D}^n(\rho))$ with symbol $\varphi^t$ (see e.g. \cite{BDF, FJ, AGS}).
\begin{definition}[Koopman semigroup \cite{ALMCM}] 
The semigroup of Koopman operators (in short, Koopman semigroup) on $\mathbb{H}^2(\mathbb{D}^n(\rho))$ is the family of linear operators $\left(U^t\right)_{t\geq 0
}$ defined by $$U^t:\mathbb{H}^2(\mathbb{D}^n(\rho))\rightarrow \mathbb{H}^2(\mathbb{D}^n(\rho)), \quad U^tf=f\circ \varphi^t.$$
\end{definition}
Under a contraction assumption on the flow $\varphi^t$, one can prove the boundedness and the strong continuity of the Koopman semigroup. In this work, we focus on the evolution of the evaluation functionals $k_z$ of the Hardy space (see \cite{CMZAM2} for the technical details), so that the above properties are not required.

\begin{definition}[Koopman generator \cite{ALMCM}, chapter 7] The Koopman generator associated with the vector field \eqref{eq:systnon} is the linear operator
\begin{equation*}
\label{eq:koopmangenerator}
L_F: \mathcal{D}(L_F) \subset \mathbb{H}^2(\mathbb{D}^n(\rho)) \rightarrow \mathbb{H}^2(\mathbb{D}^n(\rho)), \quad L_F f:=F \cdot \nabla f
\end{equation*}
with the domain \[\mathcal{D}(L_F)=\left\lbrace f\in \mathbb{H}^2(\mathbb{D}^n(\rho)): F \cdot \nabla f\in \mathbb{H}^2(\mathbb{D}^n(\rho)) \right\rbrace .\]
\end{definition}

Moreover, the expression of the Koopman generator in the basis of monomials can be obtained from the Taylor expansion 
\begin{equation}
\label{eq:Taylor}
F_l(z) = \sum_{\vert \alpha \vert\geq 1} a_{l,\alpha} \, z^\alpha = \sum_{k=1}^\infty a_{l,k} \, z^{\alpha(k)}
\end{equation}
of the vector field (with a slight abuse of notation, we will use two different conventions for the subscripts of the Taylor coefficients, i.e. $a_{l,k}=a_{l,\alpha(k)}$). It is shown in \cite{CMZAM2} that
\begin{equation}\label{eq:koop_matrix1}
\left\langle L_Fe_{k},e_{j} \right\rangle = 
\begin{cases}\sum_{l=1}^n \alpha_l(k)  \, a_{l,(\alpha(j)-\alpha(k))_l} & \textrm{if } |\alpha(j) |\geq |\alpha(k)| \\
0 &  \textrm{if } |\alpha(j) |< |\alpha(k)|.
\end{cases}
\end{equation}
with
\begin{multline*}
(\alpha(j)-\alpha(k))_l=(\alpha_1(j)-\alpha_1(k),\cdots,\alpha_l(j)-\alpha_l(k)+1, \cdots, \alpha_n(j)-\alpha_n(k)).
\end{multline*}

In particular, for monomials $e_k$ and $e_j$ of same total degree $|\alpha(j)|=|\alpha(k)|$, we have 

{\small
\begin{equation}\label{eq:total_degree_bloc}
\left\langle L_Fe_{k},e_{j} \right\rangle =
\begin{cases}\sum_{l=1}^{n}\, \alpha_l(j) \, a_{l,\alpha(l)}& \textrm{if } j=k \\
\alpha_{l}(k) \, a_{l,\alpha(r)} & \textrm{if } \alpha(j)=(\alpha_1(k),\cdots, \alpha_{l}(k)-1, \cdots, \\
&\qquad \qquad \qquad \qquad  \alpha_{r}(k)+1, \cdots,\alpha_n(k)),
\\
0 &  \textrm{otherwise}.
\end{cases}
\end{equation}}

\subsection{Stability result}

We now present an intermediate result that we will use to prove our main stability results. It is adapted from \cite{CMZAM2}, where a switched system was considered instead of \eqref{eq:systnon}.

\begin{lemma}\label{lemma:principalresult}
Consider the nonlinear system \eqref{eq:systnon} satisfying Assumptions \ref{assump4} and \ref{assump_hyperbolic} on the unit polydic. Moreover, assume that the Jacobian matrix $JF(0)$ is diagonal and there exists $\rho\in ]0, 1]$ such that $\mathbb{D}^n\left(\rho \right)$ is forward invariant with respect to the flow.
Let $\left( b_{jk}\right)_{j\geq 1, k\geq 1}$ be a double sequence of positive real numbers such that $b_{jk}b_{kj}>0$ if $\langle L_{F} e_k,  e_j \rangle \neq 0$ and such that $\sum_{k=1}^\infty b_{jk}\leq 1$, and define the double sequence $\left( Q_{jk}\right)_{j\geq 2, 1\leq k\leq j-1}$ with

\begin{equation}
\label{eq:double_sequence}
    Q_{jk} = \dfrac{\left|\left\langle L_{F} e_k,  e_j \right\rangle \right|^2  }{4\left|\Re\left(\left\langle   L_{F} e_j, e_j \right\rangle\right)\right|\left|\Re\left(\left\langle L_{F} e_k,  e_k \right\rangle\right)\right| b_{jk}b_{kj}}
\end{equation}
if $\left\langle L_{F} e_k,  e_j \right\rangle \neq 0$ and $Q_{jk}=0$ otherwise.
If the series 
\begin{equation}\label{eq:conditionconvergence_bis}
\sum_{k=1}^{+\infty}|\alpha(k)|\,\epsilon_k\,\rho^{2|\alpha(k)|}
\end{equation}
is convergent with
\begin{equation}
\label{eq:defineepsilon}
\epsilon_j > \max_{\substack{ k=1,\dots, j-1}} \epsilon_k \, Q_{jk},
\end{equation}
then the system \eqref{eq:systnon} is GAS on $\mathbb{D}^n(\rho)$. Moreover the series 
 $$V(z)=\sum_{k=1}^\infty \epsilon_k \left\vert z ^{\alpha (k)}\right\vert^2$$
 is a Lyapunov function on $ \mathbb{D}^n(\rho)$, i.e. $F(z) \cdot \nabla V(z) <0$ for all $z\in \mathbb{D}^n(\rho) \setminus\{0\}$.
\end{lemma}

The proof follows on similar lines as in \cite{CMZAM2}.

\begin{remark}
The assumption that the Jacobian matrix $JF(0)$ is diagonal can be extended to a diagonalizability condition of $JF(0)$. Indeed, if there exits $P$ such that  $J\widehat{F}(0)=P^{-1}JF(0)P$ is diagonal, a change of variables $\widehat{z}=P^{-1}z$ in $\mathbb{D}^n(\rho)$ can be chosen so that the dynamics $\dot{\widehat{z}}=\widehat{F}(\widehat{z})=P^{-1}F(P\widehat{z})$ in the new variables has a diagonal Jacobian matrix and is defined on an invariant set that is contained in $\mathbb{D}^n(\rho)$ (see the example in Section \ref{sec:ex1}). Therefore, from this point on, we will assume without loss of generality that the Jacobian matrix $JF(0)$ is diagonal.
Moreover, most of our results could be extended to upper triangular Jacobian matrices, a property which is always satisfied in $\mathbb{C}^{n \times n}$ up to a linear change of coordinates (Schur's theorem). See \cite{CMZAM2} for this general case.

\end{remark}
  
\section{Global stability criteria}
\label{sec:mainresult}

We are now in a position to present our main results. We will consider separately the case of polynomial vector fields and analytic vector fields.

\subsection{Stability criterion for polynomial vector fields}

Let us consider a dynamical system with a polynomial vector field
\begin{equation}\label{eq:syst_rescalling}
\dot{z}_l=F_l(z)=\sum_{k=1}^r a_{l,k} \, z^{\alpha(k)}, \quad l=1,\dots,n.
\end{equation}
We first define the following quantities associated with the polynomial vector field.
\begin{itemize}
\item Let $d$ be the maximal degree of the polynomials $F_l$, i.e.
\begin{equation*}
d=\max_{k \in \mathbb{N}} \left \{ |\alpha(k)| : a_{l,k} \neq 0 \textrm{ for some } l\right\} = |\alpha(r)|
\end{equation*}
\item Let $K$ be the number of nonzero terms (without counting the term containing the monomial $z_l$ in $F_l$), i.e.
\begin{equation}
\label{eq:nb_terms}
K=\sum_{l=1}^n \# \left \{ k \neq l:a_{l,k} \neq 0 \right\}
\end{equation}
where $\#$ is the cardinal of a set.
\item Let $S$ be the maximal polynomial coefficient over all components of the vector field (again discarding the terms containing the monomial $z_l$ in $F_l$), i.e.
\begin{equation*}\label{eq:cond_thrm1}
S=\max_{l=1,\cdots,n} \max_{\substack{ k=1,\cdots,r\\ k\neq l}}   \left|a_{l,k}\right|.
\end{equation*}
\item Let $R$ be the minimal real part of the diagonal entries of $JF(0)$, i.e.
\begin{equation*}\label{eq:cond_thrm12}
R=\min_{l=1,\cdots,n}   \left|\Re\left(a_{l,l}\right)\right|.
\end{equation*}
\end{itemize}

Then we have the following result.
\begin{theorem}\label{thm:poly}
Consider a dynamical system with polynomial vector field \eqref{eq:syst_rescalling} on the polydisc $\mathbb{D}^n(\mu)$, which satisfies Assumptions \ref{assump4} and \ref{assump_hyperbolic} for some $\mu>0$ large enough. Moreover, assume that the Jacobian matrix $JF(0)$ is diagonal.

 Then \eqref{eq:syst_rescalling} is GAS on $\mathbb{D}^n(\rho)$ with
 \begin{equation*}
   \rho < \begin{cases}   \dfrac{R}{KS} & \textrm{if } KS/R \geq 1 \\
    \sqrt[d-1]{\dfrac{R}{KS}} & \textrm{if } KS/R < 1
   \end{cases}
 \end{equation*}
provided that $\mu>\rho$ and $\mathbb{D}^n\left(\rho\right)$ is forward invariant with respect to the flow.
 
\end{theorem}

See Appendix \ref{sec:appendix_proof1} for the proof.

\subsection{Stability criterion for analytic vector fields}

In this section, we provide a result for dynamics with analytic vector fields, which we rewrite as
\begin{equation}\label{eq:syst_rescalling_2}
\dot{z}_l=F_l(z)=\sum_{k=1}^\infty a_{l,k} \, z^{\alpha(k)}, \quad l=1,\dots,n,
\end{equation}
under the assumption that the  Jacobian matrix $JF(0)$ is diagonal.

We first define the following quantities associated with the Taylor expansion \eqref{eq:Taylor} of the vector field.
\begin{itemize}
        \item Let $L_\mu$ be the discounted (infinite) sum of Taylor coefficients of the vector field, i.e.
        \begin{equation}\label{eq:cond_thrm2}
L_\mu=\sum_{l=1}^n  \sum_{k=1}^\infty \mu^{\vert \alpha(k) \vert} \left| a_{l,k}\right| .
\end{equation}
Note that $L_\mu$ might not be a convergent series for all $\mu$, but is always convergent for $\mu \leq 1$ under Assumption \ref{assump4}.
\item Let $R$ be the minimal real part of the diagonal entries of $JF(0)$, i.e.
\begin{equation*}
R=\min_{l=1,\cdots,n}   \left|\Re\left(a_{l,l}\right)\right|.
\end{equation*}
\end{itemize}

We have the following result.

\begin{theorem}\label{thm:gen2}
Consider a dynamical system with analytic vector field \eqref{eq:syst_rescalling_2}, which satisfies Assumptions \ref{assump4} and \ref{assump_hyperbolic}, and defined on the polydisc $\mathbb{D}^n(\mu)$ with $\mu>0$ such that $L_\mu$ is convergent. Moreover, assume that the Jacobian matrix $JF(0)$ is diagonal.

Then \eqref{eq:syst_rescalling_2} is GAS on $\mathbb{D}^n(\rho)$ with
\begin{equation}
\label{eq:criteria_thm2}
 \rho<\frac{\mu R}{L_\mu},
\end{equation}
provided that $\mathbb{D}^n\left(\rho\right)$ is forward invariant with respect to the flow.
\end{theorem}

See Appendix \ref{sec:appendix_proof2} for the proof.

\begin{remark}
If the Jacobian matrix is not diagonal(izable), the above result can be extended to the case of an upper triangular Jacobian matrix with additional diagonal dominance conditions
\begin{equation*}
\label{eq:diagonal_dominant}
   \left|a_{q,r}\right|^2 < \frac{1}{D^2} \left|\Re(a_{q,q})\right| \left| \Re(a_{r,r}) \right|, \quad 1 \leq q < r \leq n
\end{equation*}
and
\begin{equation*}
\label{eq:diagonal_dominant2}
    \left|a_{q,r}\right| < \frac{1}{D} \left|\Re(a_{q,q})\right|, \quad 1 \leq q < r \leq n
\end{equation*}
where $D$ is the number of upper off-diagonal nonzero entries of $JF(0)$.
    See the proof of Corollary 3.9 in \cite{CMZAM2} for more details.
\end{remark}

\section{Examples}
\label{sec:illustra}

In this section, we estimate the region of attraction of equilibria by using our stability criteria. We consider examples inspired by \cite{MCCDFSD}, where the authors provide some guidelines to construct vector fields that generate holomorphic flows on the bidisc $\mathbb{D}^2$.
 
\subsection{Polynomial vector field}
\label{sec:ex1}

Consider the vector field 
\begin{equation}\label{eq:example_1}
F(z_1,z_2)=\begin{cases} a \left(z_1-\frac{1}{ac}z_2\right)\\
a \left(z_2-\frac{1}{ac}z_1+bz_1^2\right),\end{cases}
\end{equation}
where $a=-1/4$, $c=8$ and $b=-1/50$.
The dynamics admit the equilibria $(0,0)$ and $\left( -75, 150\right)$ so that $(0,0)$ is the unique equilibrium point on the polydisc $\mathbb{D}^n\left(\mu\right)$ with $\mu<75$. 
The Jacobian matrix $JF(0)$ has negative eigenvalues  \mbox{$a-1/c=-3/8$} and $a+1/c=-1/8$, and is diagonalizable by the matrix $P= \begin{pmatrix} 1 & -1\\
1 & 1\end{pmatrix}$.
Using the change of coordinates $\widehat{z}=P^{-1}z$, we have
\begin{equation}\label{eq:new_example_1}
\widehat{F}(\widehat{z}_1,\widehat{z}_2)
=\begin{cases}   (a-\frac{1}{c})\widehat{z}_1 +\frac{a^2b}{2} \left( \widehat{z}_1^2-2\widehat{z}_1\widehat{z}_2+\widehat{z}_2^2 \right)   \\
(a+\frac{1}{c})\widehat{z}_2 +\frac{a^2b}{2} \left( \widehat{z}_1^2-2\widehat{z}_1\widehat{z}_2+\widehat{z}_2^2 \right).\end{cases}
\end{equation}
The dynamics in the new variables is forward invariant in the polydisc $\mathbb{D}^n\left(\widehat{\rho}\right)$, if we assume that $\widehat{\rho}=\rho\in ]1,\mu[$  since 
 \begin{itemize}
 \item {\small $|\widehat{z}_1|=\rho \Rightarrow \Re\left( \bar{\widehat{z}}_1 \widehat{F}_1(\widehat{z}) \right)=-3\rho^2/8-\rho^2\Re\left( \widehat{z}_1 \right)/1600+\rho^2\Re\left( \widehat{z}_2 \right)/800-\Re\left( \bar{\widehat{z}}_1 \widehat{z}_2^2\right)/1600<0$} as $\rho<75$ and
 \item  {\small $|\widehat{z}_2|=\rho \Rightarrow \Re\left( \bar{\widehat{z}}_2 \widehat{F}_2(\widehat{z}) \right)=-\rho^2/8-\rho^2\Re\left( \widehat{z}_2 \right)/1600+\rho^2\Re\left( \widehat{z}_1 \right)/800-\Re\left( \bar{\widehat{z}}_2 \widehat{z}_1^2\right)/1600<0$ as $\rho<75$.}
 \end{itemize}
For the vector field $\widehat{F}$, we compute $\widehat{d}=2$, $\widehat{K}=6$, $\widehat{S}=a^2|b|=1/800$ and $\widehat{R}=| a+1/c|=1/8$, so that
$\widehat{K}\widehat{S}/\widehat{R}=3/50<1.$
Hence, it follows from Theorem \ref{thm:poly} that the nonlinear system \eqref{eq:new_example_1} is GAS on $\mathbb{D}^2(\widehat{\rho})$ with $\widehat{\rho}<50/3$. Finally, this implies that \eqref{eq:example_1} is GAS on $P(\mathbb{D}^2(\widehat{\rho})) \supset \mathbb{D}^2(\widehat{\rho})=\mathbb{D}^2(\rho)$ since $\|P\|_\infty=2>1$ and with $\rho =\widehat{\rho}<50/3$.

\subsection{Analytic vector field}

Consider the vector field 

\begin{equation}\label{eq:example_2}
F(z_1,z_2)=\begin{cases} a \left(z_1-\dfrac{2z_2^2}{c-z_2}\right)\\
a \left(z_2-\dfrac{bz_1^2}{(d-z_1)^2}\right),\end{cases}
\end{equation}
where  $a=-1$, $b=4$, $c=30$ and $d=20$. The origin $(0,0)$ is the unique equilibrium point on the polydisc  $\mathbb{D}^n\left(\mu\right)$ with $\mu=10$. If we assume that  $\rho\in ]1,\mu[$, $\mathbb{D}^n\left(\rho\right)$ is invariant with respect to the flow since
 \begin{itemize}
 \item {\small $|z_1|=\rho \Rightarrow \Re\left( \bar z_1 F_1(z) \right)=-\rho^2+2\Re\left(  \dfrac{\bar z_1z_2^2}{30-z_2} \right)<0$}
 since {\small \mbox{$1> \dfrac{2\rho}{30-\rho}$}} and it follows that
{\small \[ \rho^2>\dfrac{2\rho^3}{30-\rho}>\dfrac{2\rho^3}{\left|30-|z_2|\right|}>2\left|  \dfrac{\bar z_1z_2^2}{30-z_2} \right|\geq 2\left| \Re\left(  \dfrac{\bar z_1z_2^2}{30-z_2}\right)\right| \]}
 \item {\small $|z_2|=\rho \Rightarrow \Re\left( \bar z_2 F_2(z) \right)=-\rho^2+4\Re\left( \dfrac{ \bar z_2 z_1^2}{(20-z_1)^2} \right)<0$}
  since {\small $1> \dfrac{4}{(20-\rho)^2}$} and it follows that
 {\small\begin{multline*}
      \rho^2>\dfrac{4\rho^2}{|20-\rho|^2}>\dfrac{4\rho^2}{\left|20-|z_1|\right|^2}>\left|  \dfrac{4\bar z_2z_1}{(20-z_1)^2} \right| \geq \left| \Re\left(  \dfrac{4\bar z_2z_1}{(20-z_1)^2}\right)\right| .\end{multline*}}
 \end{itemize}
As {\small\[ F_1(z)=-z_1+2\sum_{k=0}^\infty \dfrac{z_2^{k+2}}{30^{k+1}}  \text{ and }  F_2(z)=-z_2+4\sum_{k=0}^\infty \dfrac{(k+1)z_1^{k+2}}{20^{k+2}},\] }
we  obtain 
{\small\[
L_\mu=\left( 10+2 \sum_{k=0}^\infty \dfrac{10^{k+2}}{30^{k+1}} \right)+\left( 10+4 \sum_{k=0}^\infty \dfrac{(k+1)10^{k+2}}{20^{k+2}} \right)=73/3,
\]}
and  $R=1$. 
Hence, it follows from Theorem \ref{thm:gen2} that the nonlinear system \eqref{eq:example_2} is GAS on $\mathbb{D}^2(\rho)$ with $\rho<30/73$.

\section{Conclusions and future work}
\label{sec:concl}

We have obtained new sufficient conditions for global stability of nonlinear equilibria by leveraging the Koopman operator framework in the Hardy space of the polydisc. In particular, stability criteria were proposed, which provide an approximation of the region of attraction in the case of polynomial vector fields and more general analytic vector fields. 
These criteria are systematic in that they can be directly verified with the Taylor expansion coefficients of the vector field, so that they could be easily implemented in a toolbox.

We envision several perspectives for future work.
Our Koopman operator based techniques could be applied to other types of dynamical systems (e.g. limit cycles dynamics, general attractors). Moreover, our criteria seem to be conservative in some cases, so that they could be adapted to yield stability results in larger polydiscs. More importantly, the relevance and possible extension of our stability results to $\mathbb{C}^n$ could be investigated.


\appendix
\label{sec:appendix}


\section{Proof of Theorem \ref{thm:poly}}
\label{sec:appendix_proof1}

The following proof is inspired by the proof of Corollary 3.8 in \cite{CMZAM2}.

Let us consider the change of variable $z'=z/\mu$ which yields a rescaled dynamics on the unit polydisc $\mathbb{D}^n $ with the vector field
\begin{equation}\label{eq_rescal_newdynamic}
F'_l(z')=\sum_{k=1}^{r} \mu^{|\alpha(k) |-1}a_{l,k} z'^{\alpha(k)}.
\end{equation}
The Jacobian  $JF'(0)$ is also  diagonal and $K'=K$  (see \eqref{eq:nb_terms}).
In the new coordinates, the inner products  \eqref{eq:koop_matrix1} and \eqref{eq:total_degree_bloc} are given by
{\small
\begin{equation}\label{eq:koop_matrix1_new}
\left\langle L_{F'}e_{k}',e_{j}' \right\rangle = \begin{cases} \mu^{|\alpha(j) |- |\alpha(k) |} \sum_{l=1}^n \alpha_l(k)  \, a_{l,(\alpha(j)-\alpha(k))_l} & \text{if } |\alpha(j) |> |\alpha(k)| \\ 
\sum_{l=1}^n \alpha_l(k)  \, a_{l,l} & \text{if } j= k \\
0 & \text{otherwise}.
\end{cases}
\end{equation}}

Our result is proved through Lemma \ref{lemma:principalresult} with the sequence
{\small \begin{equation}\label{eq:polynomial}
\begin{cases}b_{jj}=(1-\xi)\\ b_{jk}=\dfrac{\xi}{2 K} & \textrm{if } j\neq k \textrm{ with }\left\langle  L_{F'}e_k',e_j'\right\rangle \neq 0 \, \textrm{ or } \left\langle  L_{F'}e_j',e_k'\right\rangle \neq 0 \\
 b_{jk}=0, & \textrm{if } j\neq k \textrm{ with } \left\langle  L_{F'}e_k',e_j'\right\rangle = 0\, \textrm{ or } \left\langle L_{F'}e_j',e_k'\right\rangle = 0 , \end{cases}
\end{equation}}
with $\xi\in]0,1[$. It is clear from \eqref{eq:koop_matrix1} that, for a fixed $j$ and for all $k\in \mathbb{N}\setminus\{j\}$, there are at most $K$ nonzero values $\langle  L_{F'}e_k',e_j'\rangle$ and at most $K$ nonzero values $\langle  L_{F'}e_j',e_k'\rangle$, so that the sequence \eqref{eq:polynomial} satisfies $\sum_{k=1}^\infty b_{jk} \leq 1$. The elements $Q_{jk}$ of the double sequence \eqref{eq:double_sequence} are given by
\begin{equation}\label{eq:maxim2_poly}
Q_{jk}= 
     \dfrac{K^2 \left| \left\langle  L_{F'}e_k',e_j'\right\rangle \right|^2}{\xi^2 \, \left|\Re\left(\left\langle L_{F'}e_j', e_j' \right\rangle\right)\right|\left|\Re\left(\left\langle L_{F'} e_k', e_k' \right\rangle\right)\right|} \quad j<k.
\end{equation}
Moreover, using \eqref{eq:koop_matrix1_new}, we obtain the inequalities
\begin{eqnarray*}
\left| \left\langle  L_{F'}e_k',e_j'\right\rangle \right| & \leq & \mu^{|\alpha(j) |- |\alpha(k) |}
\sum_{l=1}^n \alpha_l(k)  \left| a_{l,(\alpha(j)-\alpha(k))_l}\right|\\
 & \leq & S \mu^{|\alpha(j) |- |\alpha(k) |}
\left| \alpha(k)\right|
\end{eqnarray*}
and 
\begin{equation*}
\left|\Re\left(\left\langle L_{F'}e_j', e_j' \right\rangle\right)\right| =  \sum_{l=1}^{n}\, \alpha_l(j) \left|\Re \left(a_{l,\alpha(l)}\right) \right| \\
\geq  R \left| \alpha (j)\right|.
\end{equation*}
It follows from the above inequalities and from \eqref{eq:maxim2_poly} that
{\small 
\begin{equation*}
Q_{jk} \leq   \dfrac{ K^2 S^2 \mu^{2(|\alpha(j) |- |\alpha(k) |)}
\left| \alpha(k)\right|^2}{\xi^2 R^2 \left| \alpha (j)\right| \left| \alpha (k)\right|}
\leq   \dfrac{ K^2 S^2 }{\xi^2 R^2 } \mu^{2 \left(|\alpha(j) |- |\alpha(k) |\right)}
\end{equation*}}
where we used $\left| \alpha (j)\right| \geq \left| \alpha (k)\right|$.

If $KS/R \geq 1$, we set $\mu=1$. In this case, we have
\[ Q_{jk}\leq \dfrac{K^2S^2}{\xi^2R^2} \stackrel{\text{def}}{=}Q \]
for some $\xi\in]0,1[$. It follows that \eqref{eq:defineepsilon} is satisfied with
\begin{equation}
\label{eq:epsilon_diag1}
\epsilon_j \sim \max_{k\in \mathcal{K}_j} \left\{\epsilon_k\, Q \right\}
\end{equation}
with $\mathcal{K}_j=\{k \in\{1,\dots,j-1\} : \left\langle L_{F'} e_k',  e_j' \right\rangle \neq 0  \}$. Hence, the sequence \eqref{eq:epsilon_diag1} yields $\epsilon_j=\mathcal{O}(Q^{|\alpha(j)|})$ for $j > 1$. It follows that \eqref{eq:conditionconvergence_bis} is convergent with a radius $\rho< 1/\sqrt{Q}$, or equivalently $\rho<R/(KS)$ for some $\xi\in]0,1[$ large enough. Finally Lemma \ref{lemma:principalresult} implies that the dynamics \eqref{eq:syst_rescalling} is GAS on $\mathbb{D}^n(\rho)$.

If $KS/R<1$, we can choose $\mu>1$. In this case, we have
\[\mu^{2 \left(|\alpha(j) |- |\alpha(k) |\right)} \leq \mu^{2 \left(d- |\alpha(1) |\right)}=\mu^{2  \left(d- 1\right) } \]
and therefore
\[ Q_{jk}\leq \dfrac{K^2S^2\mu^{2\left(d- 1\right)}}{\xi^2R^2} \stackrel{\text{def}}{=}Q <1\]
for some $\xi\in]0,1[$ and with $1<\mu<\sqrt[d-1]{R/(KS)}$.
It follows that \eqref{eq:defineepsilon} is satisfied with
\begin{equation*}
\epsilon_j \sim \max_{k\in \mathcal{K}_j} \left\{\epsilon_k\, Q \right\} = 1
\end{equation*}
for $j \geq 1$. Then, \eqref{eq:conditionconvergence_bis} is convergent with a radius $\rho'< 1$ and Lemma \ref{lemma:principalresult} implies that the new dynamics $\dot{z}'=F'(z)$ is GAS on $\mathbb{D}^n(\rho')$ (note that the invariance of the new dynamics on $\mathbb{D}^n(\rho')$ directly follows from the invariance of the original dynamics on $\mathbb{D}^n(\rho)$). Hence, the orifinal dynamics \eqref{eq:syst_rescalling} is GAS on $\mathbb{D}^n(\rho)$, with $\rho=\mu \rho' <  \sqrt[d-1]{R/(KS)}$.


\section{Proof of Theorem \ref{thm:gen2}}
\label{sec:appendix_proof2}

The following proof is inspired by the proof of Corollary 3.9 in \cite{CMZAM2}.

Let us consider the change of variable $z'=z/\mu$ which yields a rescaled dynamics on the unit polydisc $\mathbb{D}^n $ with the vector field $F'(z')$ (see \eqref{eq_rescal_newdynamic} in the previous proof). In this case, the Jacobian matrix $JF'(0)$ is also diagonal.

Our result is proved through Lemma \ref{lemma:principalresult} with the sequence
{\small \begin{equation*}
\begin{cases}b_{jj}=(1-\kappa)\\ b_{jk}=0 & \textrm{if } j\neq k \textrm{ with }  |\alpha(j)|=|\alpha(k)|  \textrm{ and  } \left\langle  L_{F'}e_k',e_j'\right\rangle \neq 0\\
& \qquad \qquad \qquad \qquad \qquad \qquad \,\,\, \textrm{ or if} \left\langle  L_{F'}e_j',e_k'\right\rangle = 0\\
b_{jk}=\dfrac{\kappa}{2} \dfrac{\left| \left\langle L_{F'} e_k',e_j'\right\rangle \right|}{\sum_{l=1}^{\infty}\left| \left\langle  L_{F'}e_l',e_j'\right\rangle \right|} & \mbox{if}\,  |\alpha(k)|< |\alpha(j)|\\
 b_{jk}=\dfrac{\kappa}{2} \dfrac{\left| \left\langle  L_{F'} e_j', e_k'\right\rangle \right|}{\sum_{l=1}^{\infty}\left| \left\langle L_{F'}e_j',e_l' \right\rangle \right|} & \mbox{if}\,   |\alpha(k)|> |\alpha(j)| \end{cases}
\end{equation*}}
with  $\kappa\in ]0,1[$.  The sequence $b_{jk}$ satisfies 
\begin{multline*}
\sum_{k=1}^\infty b_{jk} <(1-\kappa) + \dfrac{\kappa}{2}    \dfrac{\sum_{k=1}^j \left| \left\langle L_{F'} e_k',e_j'\right\rangle \right|}{\sum_{l=1}^{\infty}\left| \left\langle  L_{F'}e_l',e_j'\right\rangle \right|}+ \dfrac{\kappa}{2}    \dfrac{\sum_{k=j+1}^\infty \left| \left\langle L_{F'} e_j',e_k'\right\rangle \right|}{\sum_{l=1}^{\infty}\left| \left\langle  L_{F'}e_j',e_l'\right\rangle \right|} < 1.
\end{multline*}
The elements $Q_{jk}$ of the double sequence \eqref{eq:double_sequence} are given by
{\small
\begin{equation}\label{eq:maxim2_analytic}
Q_{jk}= \begin{cases}
\dfrac{\sum_{l=1}^{\infty}\left| \left\langle  L_{F'}e_l',e_j'\right\rangle \right| \sum_{l=1}^{\infty}\left| \left\langle  L_{F'}e_k',e_l'\right\rangle \right|}{\kappa^2\left|\Re\left(\left\langle  L_{F'} e_j', e_j' \right\rangle\right)\right|\left|\Re\left(\left\langle   L_{F'} e_k',  e_k' \right\rangle\right)\right|} & \textrm{if } |\alpha(k)|\neq |\alpha(j)| \\ 
&  \textrm{ and }\left\langle L_{F'} e_k',  e_j' \right\rangle\neq 0 \\
0 & \textrm{otherwise.}
\end{cases}
\end{equation}}
We note that $\sum_{l=1}^{\infty}\left| \left\langle L_{F'}e_l', e_j' \right\rangle \right| $ and $ \sum_{l=1}^{\infty}\left| \left\langle  L_{F'}e_k',e_l'\right\rangle \right|$ are finite according to the assumptions. It is easy to see that $Q_{jk}>1$ for $|\alpha(j)| > |\alpha(k)|$. 

Moreover, with \eqref{eq:koop_matrix1},  \eqref{eq:total_degree_bloc} and \eqref{eq:cond_thrm2}, we obtain
{\small
\begin{eqnarray*}
\sum_{l=1}^\infty \left| \left\langle  L_{F'}e_l',e_j'\right\rangle \right| & \leq & \sum_{l=1}^\infty
\sum_{p=1}^n \alpha_p(l)  \left| a_{p,(\alpha(j)-\alpha(l))_p}'\right|\\
 & \leq & 
\sum_{l=1}^\infty \left| \alpha(l)\right|  \sum_{p=1}^n   \left| a_{p,(\alpha(j)-\alpha(l))_p}'\right|\\
& \leq & \left| \alpha(j)\right|  
\sum_{l=1}^\infty  \sum_{p=1}^n   \left| a_{p,(\alpha(j)-\alpha(l))_p}'\right|\\
& = & \left| \alpha(j)\right|  
\sum_{l=1}^\infty \mu^{|\alpha(j)|-|\alpha(l)|} \sum_{p=1}^n   \left| a_{p,(\alpha(j)-\alpha(l))_p}\right|\\
 & \leq & L_\mu
\left| \alpha(j)\right|,
\end{eqnarray*}}
and
{\small\begin{eqnarray*}
\sum_{l=1}^\infty \left| \left\langle  L_{F'}e_k',e_l'\right\rangle \right| & \leq & \sum_{l=1}^\infty
\sum_{p=1}^n \alpha_p(k)  \left| a_{p,(\alpha(l)-\alpha(k))_p}'\right|\\
& \leq & \left| \alpha(k)\right|  
\sum_{l=1}^\infty  \sum_{p=1}^n   \left| a_{p,(\alpha(l)-\alpha(k))_p}'\right|\\
  & = & \left| \alpha(k)\right|  
  \sum_{l=1}^\infty \mu^{|\alpha(l)|-|\alpha(k)|}  \sum_{p=1}^n \left| a_{p,(\alpha(l)-\alpha(k))_p}\right|\\
 & \leq & L_\mu
\left| \alpha(k)\right|.
\end{eqnarray*}}


It follows from the above inequalities and from \eqref{eq:maxim2_analytic} that
\[  Q_{jk}\leq \dfrac{L_\mu^2\left| \alpha (j)\right| \left| \alpha (k)\right|}{\kappa^2 R^2\left| \alpha (j)\right| \left| \alpha (k)\right|}= \dfrac{ L_\mu^2}{\kappa^2 R^2}\stackrel{\text{def}}{=}Q\]
so that \eqref{eq:defineepsilon} is satisfied with
\begin{equation}
\label{eq:epsilon_diag}
\epsilon_j \sim \max_{k\in \mathcal{K}_j} \left\{\epsilon_k \, Q \right\}
\end{equation}
 with 
$$\mathcal{K}_j=\{k \in{1,\dots,j-1} : \left\langle L_{F'} e_k',  e_j' \right\rangle \neq 0 \textrm{ for }  |\alpha(k)| < |\alpha(j)|\}.$$
Hence, the sequence \eqref{eq:epsilon_diag} yields $\epsilon_j=\mathcal{O}(Q^{|\alpha(j)|})$. It follows that \eqref{eq:conditionconvergence_bis} is convergent with a radius $\rho'<1/\sqrt{Q}$ or equivalently $\rho'<R/L_\mu$ with $\kappa \in ]0,1[$ large enough. Then Lemma \ref{lemma:principalresult} implies that the dynamics $\dot{z}'=F'(z)$  is GAS on $\mathbb{D}^n(\rho')$ (note that the invariance of the new dynamics on $\mathbb{D}^n(\rho')$ directly follows from the invariance of the original dynamics on $\mathbb{D}^n(\rho)$). Hence, the original dynamics \eqref{eq:syst_rescalling_2} is GAS on $\mathbb{D}^n(\rho)$, with $\rho=\mu \rho' <\mu R/L_\mu$.



\end{document}